\documentclass[9pt]{article}
\usepackage{amssymb}
\usepackage{amsmath}

\newtheorem{theo}{Theorem}
\newtheorem{lem}{Lemma}
\newtheorem{prop}{Proposition}
\newtheorem{cor}{Corollary}
\newtheorem{rem}{Remark}

\newtheorem{defi}{Definition}

\begin{document}

\title{Two-scale homogenization of piezoelectric perforated structures}
\date{}
\author{Houari Mechkour\thanks{%
Centre de Mathématiques Appliquées (UMR 7641)
\'Ecole Polytechnique, 91128 Palaiseau, France. (mechkour@cmap.polytechnique.fr).}}
\maketitle

\begin{abstract}
We are interested in the homogenization of  elastic-electric coupling equation, with rapidly oscillating coefficients, 
in periodically perforated piezoelectric body. We justify the two first terms in the usual asymptotic development of the problem solution. For the main convergence results of this paper, we use the notion of {\it two-scale convergence}. A two-scale homogenized 
system is obtained as the limit of the periodic problem. While in the static limit the method provides homogenized 
electroelastic coefficients whicht coincide with those deduced from other homogenization techniques (asymptotic 
homogenization \cite{Cast}, $\Gamma$-convergence \cite{Tel1}).

\noindent \textbf{Key words. }Homogenization; Piezoelectricity; Perforations
\textbf{\ }\newline
%\noindent \emph{Mathematics Subject Classification: }90C34, 49J53, 90C25,90C31.
\end{abstract}
\section{Introduction}
Composites and perforated (lattice) materials are widely used in many practical applications, such as aircraft, civil 
engineering, electrotechnics, and many others. These materials are with a large number of heterogeneities (inclusions 
or holes), and in strong contrast to continum materials, their behavior is definitively influenced by micromechanical 
events.\\

The first goal of this work we study the homogenization of the equation of the elastic-electric coupling with rapidly 
oscillating coefficients in a periodically perforated domain. The homogenized of this problem for a fixed domain has 
already been studied, by the author (Feng and Wu.\cite{Feng}, Castillero {\it and all}.\cite{Cast}, Ruan {\it and all}.
\cite{Ruan}). But in this work we give new convergence results concerning the same model by using homogenization 
technique of {\it `` two-scale convergence ''}, which permits us to conclude the limit problem, the approximation of final state is altrered by a constant named as the volum fraction which depends on the proportion of material in the perforated domain and is equal to 1 when there are no holes.\\

The second goal of this paper is to apply the technique of a formal asymptotic homogenization, to determine the 
effective elastic, piezoelectric and dielectric modulous of periodic medium. The final formulae for the effective 
parameters are given in a relatively simpler closed form.\\

The third goal of this paper is to also establish a corrector-type theorem, which permits to replace the sequence by its ``two-scale'' limit using the result on the strong convergence, and permits to justify the two first terms in the usual asymptotic expansion of the solution. In the last section we treat of the enregy aspect of our problem.
\section{Homogenization problem }
Throughout this paper $L^2(\Omega)$ in the Sobolev space of real-valued functions that are measurable and square summable in $\Omega$ with respect of the Lebesgue measure. We denote by $C^{\infty}_{\sharp}(Y)$ the space of infinitely differentiable functions in $\mathbb{R}^3$ that are periodic of $Y$. Then, $L^2_{\sharp}(Y)$ (respectively, $H^1_{\sharp}(Y)$) is the completion for the norm of $L^2(Y)$ (respectively, $H^1(Y)$) of $C^{\infty}_{\sharp}(Y)$.
\subsection{Geometric of the medium}
Let $\Omega \subset \mathbb{R}^3$ be a bounded three dimensional domain with the boundary $\Gamma = \partial \Omega$. We denote
 $x$ the macroscopic variable and by $y=\frac{x}{\varepsilon}$ the microscopic variable. Let us define 
$\Omega_{\varepsilon}$ of periodically perforated subdomains of a bounded open set $\Omega$. The period of 
$\Omega_{\varepsilon}$ is $\varepsilon Y^*$, where $Y^*$ is a subset of the unit cube $Y=(0,1)^3$, which represents 
the solid or material domain, $S^*$ obtainded by $Y$-periodicity from $Y^*$, is a smooth connected (the material is in 
one piece) open set in $\mathbb{R}^3$. Denoting by $\chi(y)$ the caracteristic function of $S^*$ (Y-periodic), in 
$\Omega_{\varepsilon}$ well be defined analytically by
$$
\Omega_\varepsilon = \big \{~x \in \Omega, ~~\chi(\frac{x}{\varepsilon}) = 1 ~\big \}
$$
\subsection{Model problem}

 We adopt the convention of Einstein for the summation of repeated indices, we use Latin indices, understood from 1 to 3, we note by ${\bf u}^{\varepsilon}$ the fields of displacement in elastic, and by $\varphi^{\varepsilon}$ of electric potentiel. The equations of equilibrium and Gauss's law of electrostatics in the absence of free charges, written as
\begin{eqnarray}
\left\{
\begin{array}{cll}
-{\bf div} \hspace{2mm}{\bf \sigma}^\varepsilon({\bf u}^\varepsilon,\varphi^\varepsilon)& = &{\bf f} ~ \mbox{in} ~ \Omega_{\varepsilon},\\
- {\bf div}\hspace{2mm}{\bf D}^\varepsilon({\bf u}^\varepsilon,\varphi^\varepsilon)& =& 0 ~\mbox{in} ~ \Omega_{\varepsilon},
\end{array}
\right.
\label{Pb1}
\end{eqnarray}
we complete the boundary conditions,
\begin{equation}
\left\{
\begin{array}{rl}
({\bf u}^{\varepsilon},\varphi^{\varepsilon})=({\bf 0},0) & \mbox{on} ~\partial\Omega, \\
{\bf \sigma}^\varepsilon({\bf u}^\varepsilon,\varphi^\varepsilon).n^\varepsilon = 0 & \mbox{on the boundary of holes} ~\partial\Omega_{\varepsilon}-\partial\Omega,\\
{\bf D}^\varepsilon({\bf u}^\varepsilon,\varphi^\varepsilon).n^\varepsilon = 0 & \mbox{on the boundary of holes} ~ \partial\Omega_{\varepsilon}-\partial\Omega,
\end{array}
\right.
\label{Pb2}
\end{equation}
where ${\bf f} \in {\bf L}^2(\Omega_{\varepsilon})$. The second-order stress tensor $\sigma^\varepsilon=(\sigma^\varepsilon_{ij})$, and the electric displacement vector ${\bf D}^\varepsilon = (D^\varepsilon_{i})$, are linearly related to the second-order strain tensor $s_{kl}({\bf u})=\frac{1}{2}(\partial_k {\bf u}_l+\partial_l {\bf u}_k)$ and the electric field vector $\partial_k \varphi^\varepsilon$ by the constitutive law
\begin{equation}
\left\{
\begin{array}{rl}
{\bf \sigma}^\varepsilon_{ij}({\bf u}^\varepsilon,\varphi^\varepsilon) = c^\varepsilon_{ijkl}s_{kl}({\bf u}^\varepsilon)+e^\varepsilon_{kij}\partial_k \varphi^\varepsilon & \mbox{in} ~\Omega_{\varepsilon}, \\
D^\varepsilon_i({\bf u}^\varepsilon,\varphi^\varepsilon) = -e^\varepsilon_{ikl}s_{kl}({\bf u}^\varepsilon)+d^\varepsilon_{ij}\partial_j \varphi^\varepsilon & \mbox{in} ~\Omega_{\varepsilon}.
\end{array}
\right.
\label{Pb3}
\end{equation}
$$
1 \leq i,j,k,l \leq 3,
$$
where $({\bf div}~\sigma^\varepsilon)^j = \partial_i \sigma^\varepsilon_{ij},~{\bf div}~{\bf D}^\varepsilon = \partial_i D^\varepsilon_i,\hspace{1mm} \partial_i = \frac{\partial}{\partial x_i}, \hspace{1mm}x = (x_i) \in \Omega$.
And the material proprieties are given by the fourth-order stiffness tensor $c^\varepsilon_{ijkl}$ measured at constant electric field, the elastic coefficients satisfy the following symmetries and ellipticity uniformily in $\varepsilon$,
\begin{equation}
\left\{
\begin{array}{ll}
c^\varepsilon_{ijkl}(x)=c_{ijkl}(x,\frac{x}{\varepsilon}),\\
c^\varepsilon_{ijkl} = c^\varepsilon_{jikl} = c^\varepsilon_{klij}=c^\varepsilon_{ijlk},\\
c_{ijkl}(x,y) \in L^{\infty}(\Omega;C_{\sharp}(Y)), \\
\exists \alpha_c \neq \alpha_c(\varepsilon) >0 : c^\varepsilon_{ijkl}X_{ij}X_{kl} \geq \alpha_c X_{ij}X_{ij}, ~~\forall X_{ij}=X_{ji} \in \mathbb{R}. 
\end{array}
\right.
\label{cond1}
\end{equation}
The third-order piezoelectric tensor $e^\varepsilon_{ijk}$ (the coupled tensor), verify the following symmetry,
\begin{equation}
\left\{
\begin{array}{ll}
e^\varepsilon_{ijk}= e_{ijk}(x,\frac{x}{\varepsilon}),\\
e^\varepsilon_{ijk} = e^\varepsilon_{ikj},\\
e_{ijk}(x,y) \in L^{\infty}(\Omega;C_{\sharp}(Y)). \\
\end{array}
\right.
\label{cond2}
\end{equation}
The second-order electric tensor $d^\varepsilon_{ij}$ (dielectric permittivity), measured at constant strain, verify the conditions of symmetric and ellipticity uniformily by $\varepsilon$,
\begin{equation}
\left\{
\begin{array}{ll}
d^\varepsilon_{ij}=d_{ij}(x,\frac{x}{\varepsilon}),\\
d^\varepsilon_{ij}=d^\varepsilon_{ji},\\
d_{ij}(x,y) \in L^{\infty}(\Omega;C_{\sharp}(Y)), \\
\exists \alpha_d\neq \alpha_d(\varepsilon) > 0 : d^\varepsilon_{ij} X_i X_j \geq \alpha_d X_i X_i,  \hspace{2mm}\forall X_i \in \mathbb{R}.
\end{array}
\right.
\label{cond3}
\end{equation}
\subsection{Variational problem}
Introducing the two Hilbert spaces 
$$
{\bf V}_{\varepsilon}(\Omega_{\varepsilon}) = \{~{\bf v} \in {\bf H}^1(\Omega_{\varepsilon}),~ {\bf v}  = {\bf 0} ~~ \mbox{on} ~~\partial\Omega  \}
$$
$$
W_{\varepsilon}(\Omega_{\varepsilon}) = \{~\psi \in H^1(\Omega_{\varepsilon}), ~\psi = 0 ~~ \mbox{on} ~~\partial\Omega   \}
$$ 
With two norms : $\parallel . \parallel_{{\bf V}_{\varepsilon}(\Omega_{\varepsilon})} =\parallel . \parallel_{{\bf H}^1(\Omega)},~~ \parallel . \parallel_{W_{\varepsilon}(\Omega_{\varepsilon})} =\parallel . \parallel_{H^1(\Omega)}$. The variational problem is defined by :
\begin{equation}
\left\{
\begin{array}{lcl}
\mbox{Find} \hspace{2mm}({\bf u}^{\varepsilon}, \varphi^{\varepsilon}) \in {\bf V}_{\varepsilon}(\Omega_{\varepsilon})\times W_{\varepsilon}(\Omega_{\varepsilon}),~\mbox{such that}\\ \\
a_{\varepsilon}(({\bf u}^{\varepsilon}, \varphi^{\varepsilon}),({\bf v},\psi))=L_{\varepsilon}({\bf v},\psi)\hspace{1cm}\forall ({\bf v},\psi)\in {\bf V}_{\varepsilon}(\Omega_{\varepsilon}) \times W_{\varepsilon}(\Omega_{\varepsilon}),
\end{array}
\right.
\label{PV1}
\end{equation}
where
\begin{eqnarray}
\left\{
\begin{array}{lcl}
\displaystyle a_{\varepsilon}(({\bf u}^{\varepsilon}, \varphi^{\varepsilon}),({\bf v},\psi))&=&\displaystyle  \int_{\Omega_{\varepsilon}}\{[c^\varepsilon_{ijkl}s_{kl}({\bf u}^\varepsilon)+e^\varepsilon_{kij}\partial_k \varphi^\varepsilon]s_{ij}({\bf v})\\\\
&\displaystyle+&[-e^\varepsilon_{ikl}s_{kl}({\bf u}^\varepsilon)+d^\varepsilon_{ij}\partial_j \varphi^\varepsilon]\partial_i \psi \}~dx\\
\displaystyle L_{\varepsilon}({\bf v},\psi) = \int_{\Omega_{\varepsilon}} f_i ~ v_i~dx
\end{array}
\right.
\label{PV2}
\end{eqnarray}
It is pointed out that under assumptions (\ref {cond1})-(\ref {cond2})-(\ref {cond3}), the variational problem 
(\ref {PV1})-(\ref {PV2}) have a unique solution $ ({\bf u}^{\varepsilon},\varphi^{\varepsilon}) \in {\bf V}_{\varepsilon}(\Omega_{\varepsilon}) \times W_{\varepsilon}(\Omega_{\varepsilon})$, 
corresponding {\it the saddle point} of this functional (see \cite{Mech2}) :
$$
({\bf v},\psi)\rightarrow \frac{1}{2}\int_{\Omega_{\varepsilon}} (c^{\varepsilon}({\bf v},{\bf v})+2e^{\varepsilon}({\bf u},\psi)-d^{\varepsilon}(\psi,\psi))~dx -\int_{\Omega_{\varepsilon}} {\bf f}~{\bf v}~dx,
$$
where
$$
\left\{
\begin{array}{lcl}
c^{\varepsilon}({\bf u},{\bf v})&=&c^\varepsilon_{ijkl}~s_{ij}({\bf u})~s_{kl}({\bf v})\\
e^{\varepsilon}({\bf u},\psi)&=&e^\varepsilon_{ikl}~s_{kl}({\bf u})~\partial_i\psi\\
d^{\varepsilon}(\psi,\psi)&=&d^{\varepsilon}_{ij}~\partial_i \psi~\partial_j \psi
\end{array}
\right.
$$ 
\subsection{A priori estimates}
 In order to prove the main convergence results of this paper we use the notion of {\it two-scale convergence} which was introduced in \cite{Ngue} and developed further in \cite{All2}. The idea of this convergence is based in first step by taking a priori estimates for displacement field and the electric potentiel. The second step we use the relatively compact property with the classical procedure of prolongation (wich is the extension by $0$ from $\Omega_{\varepsilon}$ to $\Omega$). Finally we pass the limit $\varepsilon \rightarrow 0$, in order to obtain the homogenized and the local problems in same time.\\

\begin{prop}
{\it
By using the two equivalent norms of ${\bf V}_{\varepsilon}(\Omega_{\varepsilon}),~W_{\varepsilon}(\Omega_{\varepsilon})$, for any sequence of solution $({\bf u}^{\varepsilon},\varphi^{\varepsilon})_{\varepsilon} \subset {\bf V}_{\varepsilon}(\Omega_{\varepsilon}) \times W_{\varepsilon}(\Omega_{\varepsilon})$ of variational problem (\ref {PV1})-(\ref {PV2}). Then, this solution is bounded, and we have this a priori estimate uniformly by $\varepsilon$

\begin{equation}
\parallel {\bf u}^{\varepsilon}  \parallel_{{\bf H}^1(\Omega_{\varepsilon})} + \parallel \varphi^{\varepsilon}  \parallel_{H^1(\Omega_{\varepsilon})} \leq C,
\label{estim}
\end{equation}
where $C$ is constant strictement positive and independent by $\varepsilon$.
}
\end{prop}
{\bf Proof.} \\
By choosing ${\bf v}={\bf u}^{\varepsilon}$ and $\psi = \varphi^{\varepsilon}$ in variational formulae 
(\ref {PV1})-(\ref {PV2}), and by using the Korn's and Poincar\'e's inequalities in perforated domains ({\it see } 
Oleinik {\it et al.} \cite{Olei} for the Korn's inequality and Allaire-Murat \cite{All3} for Poincar\'e's inequality), 
we see that ${\bf u}^{\varepsilon}$ and $\varphi^{\varepsilon}$ are bounded, by a constant which does not depend on 
$\varepsilon$. For other details see \cite{Mech2}.
 \section{Two-scale convergence}
We denote by $\stackrel{\sim}{.}$ the extension by zero in the holes $\Omega-\Omega_{\varepsilon}$. The sequence of solution $({\bf u}^{\varepsilon},\varphi^{\varepsilon})_{\varepsilon} \subset {\bf V}_{\varepsilon}(\Omega_{\varepsilon}) \times W_{\varepsilon}(\Omega_{\varepsilon})$ of variational problem (\ref {PV1})-(\ref {PV2}) verify (\ref {estim}) and, in this case, by adding the relatively compact property and elementary properties of two-scale convergence, imply\\
\begin{lem}
{\it
\begin{enumerate}
\item{There exists ${\bf u}(x)\in {\bf H}^1_0(\Omega)$ and $\varphi(x) \in H^1_0(\Omega)$ such that, the two sequences $(\stackrel{\sim}{{\bf u}^{\varepsilon}})_{\varepsilon},~(\stackrel{\sim}{\varphi^{\varepsilon}})_{\varepsilon}$ two-scale converge to $\chi(y){\bf u}(x),~\chi(y)\varphi(x)$, respectively.}
\item{There exists ${\bf u}_1(x,y)\in {\bf L}^2[\Omega;{\bf H}^1_{\sharp}(Y^{*})/\mathbb{R}],~\varphi_1(x,y)\in L^2[\Omega;H^1_{\sharp}(Y^{*})/\mathbb{R}]$ such that, 
$$\stackrel{\sim}{\nabla} {\bf u}^{\varepsilon} \rightarrow \chi(y)[\nabla_{x}{\bf u}(x)+\nabla_{y}{\bf u}_1(x,y)]~~ \mbox{in two-scale sense}$$
$$\stackrel{\sim}{\nabla} \varphi^{\varepsilon} \rightarrow \chi(y)[\nabla_{x}\varphi(x)+\nabla_{y}\varphi_1(x,y)]~~ \mbox{in two-scale sense}$$}
\item{ We have
$$\stackrel{\sim}{s}({\bf u}^{\varepsilon}) \rightarrow \chi(y)[s_{x}({\bf u}(x))+s_{y}({\bf u}_1(x,y))]~~ \mbox{in two-scale sense}$$}
index $x$ or $y$ means that the derivatives are with respect to the variable.
\end{enumerate}
}
\end{lem}
{\bf Proof.} For details see \cite{All2}, \cite{Ngue}, \cite{Mech1} and \cite{Mech2}.\\
\begin{cor}
{\it The sequence $(\stackrel{\sim}{{\bf u}^{\varepsilon}})_{\varepsilon>0}$ $( \mbox{resp}.~ (\stackrel{\sim}{{\varphi}^{\varepsilon}})_{\varepsilon>0})$, converge weakly to a limit $\theta {\bf u}$ $(\mbox{resp}.~\theta \varphi)$ in ${\bf L}^2(\Omega)$ $(\mbox{resp}.~ L^2(\Omega))$.}
\end{cor}
\begin{rem} 
Let $\rho \in L^2_{\sharp}(Y)$, define $\rho^{\varepsilon}(x) = \rho(\frac{x}{\varepsilon})$, and 
$(v^{\varepsilon})_{\varepsilon} \subset L^2(\Omega)$ two-scale converge to a limit $v \in L^2(\Omega \times Y)$. 
Then $(\rho^{\varepsilon}v^{\varepsilon})_{\varepsilon}$ two-scale converges to a limit $\rho v$ (see \cite{Mech2}).
\end{rem}
From last results, we can state next theorem
\begin{theo}
{\it 
The sequences $(\stackrel{\sim}{{\bf u}^{\varepsilon}})_{\varepsilon}$, $( \stackrel{\sim}{s}({\bf u}^{\varepsilon}))_{\varepsilon}$, $(\stackrel{\sim}{\varphi^{\varepsilon}})_{\varepsilon}$ and $(\stackrel{\sim}{\nabla}\varphi^{\varepsilon})_{\varepsilon}$ two-scale converge to $\chi(y){\bf u}(x)$, $\chi(y)[s_{x}({\bf u})+s_{y}({\bf u}_1)]$, $\chi(y)\varphi(x)$ and $\chi(y)[\nabla_{x} \varphi +\nabla_{y}\varphi_1]$ respectively, where $({\bf u}(x),~{\bf u}_1(x,y),~\varphi(x),~\varphi_1(x,y))$ are the unique solutions in ${\bf H}^1_0(\Omega) \times {\bf L}^2[\Omega;{\bf H}^1_{\sharp}(Y^{*})/\mathbb{R}] \times H^1_0(\Omega) \times L^2[\Omega;H^1_{\sharp}(Y^{*})/\mathbb{R}] $ of the following two-scale homogenized system; 
\begin{equation}
\left\{
\begin{array}{lclll}
\displaystyle-\frac{\partial}{\partial x_j}[\int_{Y^*}\{ c_{ijkl}(x,y)[s_{kl,x}({\bf u})+s_{kl,y}({\bf u}_1)] + e_{kij}(x,y)[\partial_{k,x}\varphi+\partial_{k,y}\varphi_1] \}dy ] =\theta f_i(x)\mbox{ in } \Omega,\\ \\ \\ 
\displaystyle-\frac{\partial}{\partial x_i}[\int_{Y^*}\{ -e_{ikl}(x,y)[s_{kl,x}({\bf u})+s_{kl,y}({\bf u}_1)] + d_{ij}(x,y)[\partial_{j,x}\varphi+\partial_{j,y}\varphi_1]\}dy ] = 0 \mbox{ in } \Omega,\\ \\ \\
\displaystyle-\frac{\partial}{\partial y_j}\{ c_{ijkl}(x,y)[s_{kl,x}({\bf u})+s_{kl,y}({\bf u}_1)] + e_{kij}(x,y)[\partial_{k,x}\varphi+\partial_{k,y}\varphi_1] \} = 0 \mbox{ in } \Omega \times Y^*,\\ \\ \\
\displaystyle-\frac{\partial}{\partial y_i}\{ -e_{ikl}(x,y)[s_{kl,x}({\bf u})+s_{kl,y}({\bf u}_1)] + d_{ij}(x,y)[\partial_{j,x}\varphi+\partial_{j,y}\varphi_1]\}= 0\mbox{ in } \Omega \times Y^*,
\end{array}
\right.
\label{pbhom1}
\end{equation}  
and we have this boundary conditions
\begin{equation}
\left\{
\begin{array}{lclll}
{\bf u}(x)& = &{\bf 0} & \mbox{on} ~ \partial \Omega, \\ \\ 
\varphi(x)& =& 0 & \mbox{on} ~\partial \Omega, \\ \\ 
\{c_{ijkl}(x,y)[s_{kl,x}({\bf u})+s_{kl,y}({\bf u}_1)] + e_{kij}(x,y)[\partial_{k,x}\varphi+\partial_{k,y}\varphi_1]\}.n_j& =& 0 & \mbox{on} ~\partial Y^*-\partial Y,\\ \\ 
\{ -e_{ikl}(x,y)[s_{kl,x}({\bf u})+s_{kl,y}({\bf u}_1)] + d_{ij}(x,y)[\partial_{j,x}\varphi+\partial_{j,y}\varphi_1]\}.n_i& =& 0 & \mbox{on} ~\partial Y^*-\partial Y.
\end{array}
\right.
\label{pbhom2}
\end{equation}
\begin{equation}
\left\{
\begin{array}{lcl}
y \rightarrow {\bf u}_1(x,y) & \mbox{is}~ Y-\mbox{periodic}, \\ \\ 
y \rightarrow \varphi_1(x,y) & \mbox{is}~ Y-\mbox{periodic}, 
\end{array}
\right.
\label{pbhom3}
\end{equation}
where $\theta$ is the volum fraction of material $($$i.e.$ $\theta=<\chi> =\int_{Y}\chi(y)~dy=\mid Y^* \mid$$)$, 
$\mid Y \mid$ denote mesure of $Y$. 
}
\end{theo}
The equations (\ref {pbhom1})-(\ref {pbhom2})-(\ref {pbhom3}) are referred to as the {\bf two-scale homogenized system}.
\paragraph{Proof.}
From the idea of G.Nguetseng \cite{Ngue}, the test functions in (\ref {PV1})-(\ref {PV2}) are chossed on the form 
$$ {\bf v}^{\varepsilon}(x) = {\bf v}(x,\frac{x}{\varepsilon})=~{\bf v}^0(x) + \varepsilon {\bf v}^1(x,\frac{x}{\varepsilon}), $$
$$ \psi^{\varepsilon}(x) = \psi(x,\frac{x}{\varepsilon})=~\psi^0(x) + \varepsilon \psi^1(x,\frac{x}{\varepsilon}), $$
where ${\bf v}^0\in {\bf C}^{\infty}_0(\Omega),~\psi^0~\in C^{\infty}_0(\Omega),~{\bf v}^1\in {\bf C}^{\infty}_0(\Omega;{\bf C}^{\infty}_{\sharp}(Y)) \mbox{ and } \psi^1~\in C^{\infty}_0(\Omega;C^{\infty}_{\sharp}(Y))$, we obtain
\begin{eqnarray*}
&\displaystyle\int_{\Omega_{\varepsilon}}\{[c^\varepsilon_{ijkl}s_{kl}({\bf u}^\varepsilon)+e^\varepsilon_{kij}\partial_k \varphi^\varepsilon][s_{ij,x}({\bf v}^0)(x)+\{s_{ij,y}({\bf v}^1)+\varepsilon s_{ij,x}({\bf v}^1)\}(x,\frac{x}{\varepsilon})]&\nonumber \\
&\displaystyle-[-e^\varepsilon_{kij}s_{kl}({\bf u}^\varepsilon)+d^\varepsilon_{ij}\partial_j \varphi^\varepsilon][\partial_{i,x} \psi^0(x)+\{\partial_{i,y} \psi^1+\varepsilon \partial_{i,x}\psi^1 \}(x,\frac{x}{\varepsilon})]\}~dx& \nonumber \\
&=\displaystyle \int_{\Omega_{\varepsilon}} f_i(x)~[v_i^0(x)+\varepsilon  v^1_i(x,\frac{x}{\varepsilon})]dx&
\end{eqnarray*}
Under the precedent hypotheses, and passing to the two-scale limit, yields
\begin{eqnarray}
&\displaystyle\int_{\Omega}\int_{Y} [c_{ijkl}(x,y)\chi(y)(s_{kl,x}({\bf u})+s_{kl,y}({\bf u}_1))+e_{kij}(x,y)\chi(y)(\partial_{k,x}\varphi+\partial_{k,y}\varphi_1)]\chi(y)&\nonumber \\
&[\displaystyle s_{kl,x}({\bf v}^0)+s_{kl,y}({\bf v}^1)]~dx~dy&\nonumber \\
&\displaystyle-\int_{\Omega}\int_{Y}[-e_{ikl}(x,y)\chi(y)(s_{kl,x}({\bf u})+s_{kl,y}({\bf u}_1))+d_{ij}(x,y)\chi(y)(\partial_{j,x}\varphi+\partial_{j,y}\varphi_1)]\chi(y)&\nonumber \\
&\displaystyle[\partial_{i,x}\psi^0+\partial_{i,y}\psi^1]~dx~dy&\nonumber \\
&\displaystyle= \int_{\Omega}\int_{Y} f_i(x) \chi(y) v_i^0(x) ~dx~dy.&
\end{eqnarray}
By definition of $\chi$, we have
\begin{eqnarray}
&\displaystyle \int_{\Omega}\int_{Y^{*}} [c_{ijkl}(x,y)(s_{kl,x}({\bf u})+s_{kl,y}({\bf u}_1))+e_{kij}(x,y)(\partial_{k,x}\varphi+\partial_{k,y}\varphi_1)]&\nonumber\\
&\displaystyle [s_{kl,x}({\bf v}^0)+s_{kl,y}({\bf v}^1)]~dx~dy&\nonumber \\
&\displaystyle +\int_{\Omega}\int_{Y^{*}} [-e_{ikl}(x,y)(s_{kl,x}({\bf u})+s_{kl,y}({\bf u}_1))+d_{ij}(x,y)(\partial_{j,x}\varphi+\partial_{j,y}\varphi_1)]&\nonumber \\
\displaystyle& [\partial_{i,x}\psi^0+\partial_{i,y}\psi^1]~dx~dy&\nonumber \\
&=\displaystyle \theta \int_{\Omega}f_i(x)v_i(x)~dx,&
\label{rttt}
\end{eqnarray}
where $\theta = \int_{Y}\chi(y)~dy$, by density of spaces from which we chose the test functions, the equation (\ref {rttt}) holds true for any ${\bf v}^0\in {\bf H}^1_0(\Omega),~\psi^0 \in H^1_0(\Omega)$, and for any  ${\bf v}^1\in {\bf L}^2[\Omega;{\bf H}^1_{\sharp}(Y^*)/\mathbb{R}]$, $\psi^1 \in L^2[\Omega;H^1_{\sharp}(Y^*)/\mathbb{R}]$,\\ 
Integrating by parts, shows that (\ref {rttt}) is variational formulation associated to the {\it two-scale homogenized system}
\begin{equation}
\left\{
\begin{array}{lcl}
-\displaystyle\partial_j \big[\int_{Y^{\star}} \{c_{ijkl}(x,y)[s_{kl,x}({\bf u})+s_{kl,y}({\bf u}_1)]+e_{kij}(x,y)[\partial_{k,x}\varphi+\partial_{k,y}\varphi_1\}~dy \big] =\theta f_i(x)\\ \\
\displaystyle-\partial_{i}\big[\int_{Y^{\star}}\{-e_{ikl}(x,y)[s_{kl,x}({\bf u})+s_{kl,y}({\bf u}_1)]+d_{ij}(x,y)[\partial_{j,x}\varphi+\partial_{j,y}\varphi_1]\}~dy \big]=0
\end{array}
\right.
\label{titi}
\end{equation}
We complete (\ref {titi}) by the boundary conditions (\ref {pbhom2})-(\ref {pbhom3}). To prove existence and uniqueness in (\ref {rttt}), 
by application of the Lax-Milgram lemma, let focus on the coercivity in 
${\bf H}^1_0(\Omega)\times {\bf L}^2[\Omega;{\bf H}^1_{\sharp}(Y^*)/\mathbb{R}] \times H^1_0(\Omega) \times L^2[\Omega;H^1_{\sharp}(Y^*)/\mathbb{R}]$ of the bilinear 
form defined by the left-hand side of (\ref {rttt}) (For a complete demonstration  see \cite{Mech2}). \\
\begin{rem}
 It is evident that the two-scale homogenized problem (\ref {pbhom1})-(\ref {pbhom2})-(\ref {pbhom3}) is a system of four equation, four unknown $({\bf u},~{\bf u}_1,~\varphi,~\varphi_1)$, each dependent on both space variables $x$ and $y$ (i.e. the macroscopic and microscopic scales) which are mixed. Although seems to be complicated, it is well-posed system of equations. Also it is clear that the two-scale homogenized problem has the same form as the original equation.\\
\end{rem}
 The object of new paragraph is to give another form of theorem which is more suitable for further physical interpretations. Indeed, we shall eliminate the microscopic variable $y$ (one doesn't want to solve the small scale structure), and decouple the two-scale homogenized problem (\ref {pbhom1})-(\ref {pbhom2})-(\ref {pbhom3}) in homogenized and cell equations. However, it is preferable, from a physical or numerical point of view ({\it see} \cite{Mech2}).
\section{Derivation of the homogenized coefficients}
Due to the linearity of the original problem, and assuming the regularity in variation of the coefficients, we take
\begin{equation}
{\bf u}_1(x,y) = s_{mh,x}({\bf u}(x)){\bf w}^{mh}(y)+\frac{\partial \varphi(x)}{\partial x_{n}} {\bf q}^{n}(y),
\label{e1}
\end{equation}
\begin{equation}
\varphi_1(x,y) = s_{mh,x}({\bf u}(x)) \varphi^{mh}(y)+\frac{\partial \varphi(x)}{\partial x_{n}} \psi^{n}(y),
\label{e2}
\end{equation}
where ${\bf w}^{mh},~\varphi^n,~{\bf q}^{mh}$ and $\psi^n$ are $Y^*$-periodic functions in $y$, independent of $x$, solutions of these two locals problems in 
$Y^*$
\begin{eqnarray}
\left\{
\begin{array}{lclllll}
\displaystyle-\frac{\partial}{\partial y_j}\Big \{ c_{ijkl}(x,y)\Big [\tau^{kl}_{mh}+ s_{kl,y}({\bf w}^{mh})\Big ]+e_{kij}(x,y) \frac{ \partial \varphi^{mh}}{\partial y_k} \Big\}&=0~\mbox{in}~Y^*, \\ \\
\displaystyle-\frac{\partial}{\partial y_i}\Big\{-e_{ikl}(x,y)\Big [\tau^{kl}_{mh}+s_{kl,y}({\bf w}^{mh})\Big ]+d_{ij}(x,y)\frac{\partial \varphi^{mh}}{\partial y_j}\Big\}&=0~\mbox{in}~Y^*,\\ \\
\displaystyle{\bf w}^{mh},~~\varphi^{mh}&Y^*-\mbox{periodics},
\end{array}
\right.
\label{pb1}
\end{eqnarray}
where
$$
\tau^{kl}_{mh}=\frac{1}{2}[\delta_{km}\delta_{lh}+\delta_{kh}\delta_{lm}]~~~1 \leq k,m,l,h \leq 3.
$$
\begin{eqnarray}
\left\{
\begin{array}{lclllll}
\displaystyle-\frac{\partial}{\partial y_j}\Big\{c_{ijkl}(x,y)s_{kl,y}({\bf q}^{n})+e_{kij}(x,y)\Big [\delta_{kn}+\frac{\partial \psi^{n}}{\partial y_k}\Big ]\Big\}&=0~\mbox{in}~Y^*,\\ \\
\displaystyle-\frac{\partial}{\partial y_i}\Big\{-e_{ikl}(x,y)s_{kl,y}({\bf q}^{n}) +d_{ij}(x,y)\Big [\delta_{jn}+\frac{\partial \psi^{n}}{\partial y_j}\Big ]\Big\}&=0~\mbox{in}~Y^*,\\ \\
\displaystyle \varphi^{n},~~\psi^{n}&Y^*-\mbox{periodics}.
\end{array}
\right.
\label{pb2}
\end{eqnarray}
However, in general a relation (\ref {e1})-(\ref {e2}) like this does not exist, if we have't linearity of problem. \\
Now substitue the expansions (\ref {e1}) and (\ref {e2}) in this equation
\\ \\
$$-\frac{\partial}{\partial y_j}\{ c_{ijkl}(x,y)[s_{kl,x}({\bf u})+s_{kl,y}({\bf u}_1)] + e_{kij}(x,y)[\partial_{k,x}\varphi+\partial_{k,y}\varphi_1] \} = 0 ~~ \mbox{in} ~\Omega \times Y^*,$$
we obtain
\begin{eqnarray}
&-&\displaystyle\frac{\partial s_{mh,x}({\bf u}) }{\partial y_j}\Big\{c_{ijkl}(x,y)\Big [\tau^{kl}_{mh}+ s_{kl,y}({\bf w}^{mh})\Big ]+e_{kij}(x,y)\frac{ \partial \varphi^{mh}}{\partial y_k} \Big\}\nonumber \\
&+&\displaystyle\frac{\partial \varphi}{\partial x_n}\Big\{c_{ijkl}(x,y)s_{kl,y}({\bf q}^{n})+e_{kij}(x,y)\Big [\delta_{kn}+\frac{\partial \psi^{n}}{\partial y_k}\Big ]\Big\} =0.
\label{kio1}
\end{eqnarray}
Calling $\tau_{mh}$ the basic of symmetric second order tensors $\tau^{kl}_{mh}=\frac{1}{2}[\delta_{km}\delta_{lh}+\delta_{kh}\delta_{lm}]$, where $\delta_{ij}$ is the Kronecker symbol.\\
Analogoulosy, we substitue the expansions (\ref {e1}) and (\ref {e2}) in this equation
$$
-\frac{\partial}{\partial y_i}\{ -e_{ikl}(x,y)[s_{kl,x}({\bf u})+s_{kl,y}({\bf u}_1)]+ d_{ij}(x,y)[\partial_{j,x}\varphi+\partial_{j,y}\varphi_1]\} = 0 ~ \mbox{in} ~~\Omega \times Y^*,
$$
we obtain
\begin{eqnarray}
&-& \displaystyle\frac{\partial s_{mh,x}({\bf u})}{\partial y_i}\Big\{-e_{ikl}(x,y)\Big [\tau^{kl}_{mh}+s_{kl,y}({\bf w}^{mh})\Big ]+d_{ij}(x,y) \frac{\partial \varphi^{mh}}{\partial y_j}\Big\}\nonumber \\
&+&\displaystyle \frac{\partial \varphi}{\partial x_{n}} \Big\{-e_{ikl}(x,y)s_{kl,y}({\bf q}^{n}) +d_{ij}(x,y)\Big [\delta_{jn}+\frac{\partial \psi^{n}}{\partial y_j}\Big ]\Big\}=0.
\label{kio2}
\end{eqnarray}
From the relation (\ref {kio1})-(\ref {kio2}), after lenghty calculations we arrive at the homogenized (effective) coefficients :

\begin{equation}
\displaystyle c_{ijmh}^H =  \Big \langle c_{ijkl}(x,y)\Big [\tau^{kl}_{mh}+s_{kl,y}({\bf w}^{mh})\Big ]+e_{kij}(x,y)\frac{\partial \varphi^{mh}}{\partial y_k} \Big\rangle,
\label{c}
\end{equation}
\begin{equation}
\displaystyle e_{nij}^H  = \Big\langle c_{ijkl}(x,y) s_{kl,y}({\bf q}^{n})+e_{kij}(x,y)\Big [\delta_{kn}+\frac{\partial \psi^{n}}{\partial y_k}\Big ] \Big\rangle, 
\label{e}
\end{equation}
\begin{equation}
\displaystyle f_{imh}^H = \Big\langle e_{ikl}(x,y)\Big [\tau^{kl}_{mh}+s_{kl,y}({\bf w}^{mh})\Big ]-d_{ij}(x,y)\frac{\partial \varphi^{mh}}{\partial y_j} \Big\rangle,
\label{f}
\end{equation}
\begin{equation}
\displaystyle d_{in}^H =\Big\langle -e_{ikl}(x,y) s_{kl,y}({\bf q}^{n}) +d_{ij}(x,y)\Big [\delta_{jn}+\frac{\partial \psi^{n}}{\partial y_j}\Big ] \Big\rangle,
\label{d}
\end{equation}
where $\langle h \rangle = \int_{Y^*} h(y)~ dy, $ the measurements on $Y^*$ of function $h$.\\
Now we give the results concerning some properties of elasticity homogenized tensor.
\begin{prop}
The coefficients of elasticity homogenized tensor $\mathcal{C^H}=(c^H_{ijkl})$ defined by (\ref{c}) satisfy:
\\ \\
a) {$c^H_{ijkl}=c^H_{klij}=c^H_{ijlk}=c^H_{jilk}, \quad \forall 1 \leq i,j,k,l \leq 3,$}
\\
b) There exists $\alpha^H_c >0$, such that for all $\xi$, symmetric tensor ($\xi_{ij}=\xi_{ji}$), 
$$c^H_{ijkl}~\xi_{ij}~\xi_{kl} \geq \alpha^H_c~\xi_{ij}~\xi_{ij} $$
\end{prop}
{\bf Proof.}\\
The part of the symmetry of these coefficients is evident
$$c^H_{ijmh}=c^H_{jimh}=c^H_{jihm}.$$
We are interesing in the proof of 
$$c^H_{ijmh}=c^H_{mhij}$$ 
Following the ideas, we transform the above expression to obtain a symmetric form. \\
We defined the tensor of second order $\Sigma$, by $\Sigma^{kl}=\frac{1}{2}(y_k \vec{e}_l+y_l \vec{e}_k)$, we define $3 \times 3$ 
$H^{kl}$ matrix by $H^{kl}=s_y(\Sigma^{kl})$, it is evident the coefficients of this matrix its defined by
$$
[H^{kl}]_{mh}=\tau^{kl}_{mh}=\frac{1}{2}[\delta_{km}\delta_{lh}+\delta_{kh}\delta_{lm}]~~~1 \leq k,l,m,h \leq 3.
$$
If we use this new notation, we can rewrite the problem (\ref {pb1}), as the form given as under
\begin{eqnarray}
\left\{
\begin{array}{lclll}
\displaystyle-\frac{\partial}{\partial y_j}\Big\{ c_{ijkl}(x,y)s_{kl,y}\Big (\Sigma^{mh}+ {\bf w}^{mh}\Big )+e_{kij}(x,y) \frac{ \partial \varphi^{mh}}{\partial y_k} \Big\}&=0~\mbox{in}~Y^*, \\ \\
\displaystyle-\frac{\partial}{\partial y_i}\Big\{-e_{ikl}(x,y)s_{kl,y}\Big (\Sigma^{mh}+ {\bf w}^{mh}\Big )+d_{ij}(x,y) \frac{\partial \varphi^{mh}}{\partial y_j}\Big\}&=0~\mbox{in}~Y^*,\\ \\
\displaystyle{\bf w}^{mh},~~\varphi^{mh}&Y^*-\mbox{periodics}.
\end{array}
\right.
\label{pba}
\end{eqnarray}
We introduce the problem functions $({\bf w}^{ij},{\bf q}^{ij})$, solutions of the problem
\begin{eqnarray}
\left\{
\begin{array}{lcllll}
\displaystyle-\frac{\partial}{\partial y_j}\Big\{ c_{kj\alpha\beta}(x,y)s_{\alpha\beta,y}\Big (\Sigma^{ij}+ {\bf w}^{ij}\Big )+e_{\alpha kj}(x,y) \frac{ \partial \varphi^{ij}}{\partial y_\alpha} \Big\}&=0~\mbox{in}~Y^*, \\ \\
\displaystyle-\frac{\partial}{\partial y_k}\Big\{-e_{k\alpha\beta}(x,y)s_{\alpha\beta,y}\Big (\Sigma^{ij}+ {\bf w}^{ij}\Big )+d_{kj}(x,y) \frac{\partial \varphi^{ij}}{\partial y_j}\Big\}&=0~\mbox{in}~Y^*,\\ \\
\displaystyle{\bf w}^{ij},~~\varphi^{ij}&Y^*-\mbox{periodics}.
\end{array}
\right.
\label{pbb}
\end{eqnarray}
The coefficient of elasticity tensor can be rewritten as
\begin{equation}
c_{ijmh}^H  = \int_{Y^*}c_{ijkl}(x,y)s_{kl,y}\Big(\Sigma^{mh}+ {\bf w}^{mh}\Big)~dy+\int_{Y^*}e_{kij}(x,y) \frac{ \partial \varphi^{mh}}{\partial y_k}~dy.
\label{Ch}
\end{equation}
The second integral of the right-hand side of precedent expression, is evaluated as follows 
\begin{eqnarray}
\int_{Y^*}e_{kij}(x,y) \frac{ \partial \varphi^{mh}}{\partial y_k}~dy&=&\int_{Y^*}e_{k\alpha\beta}(x,y) \frac{ \partial \varphi^{mh}}{\partial y_k}\delta_{\alpha i}\delta_{\beta j}~dy,\nonumber\\
&=&\frac{1}{2}\int_{Y^*}e_{k\alpha\beta}(x,y)\frac{ \partial \varphi^{mh}}{\partial y_k} ~\Big (\delta_{\alpha i}\delta_{\beta j}+\delta_{\alpha j}\delta_{\beta i}\Big )~dy,\nonumber\\
&=&-\int_{Y^*}e_{k\alpha\beta}(x,y)\frac{ \partial \varphi^{mh}}{\partial y_k} s_{\alpha \beta,y}({\bf w}^{ij})~dy\nonumber \\
&+&\int_{Y^*}e_{k\alpha\beta}(x,y)\frac{ \partial \varphi^{mh}}{\partial y_k} s_{\alpha \beta,y}(\Sigma^{ij}+{\bf w}^{ij})~dy.
\end{eqnarray}
We use the variationnal formulation of the first equation of problem (\ref {pba}), and taking the test function 
${\bf v}={\bf w}^{ij}$, we obtain 
\begin{equation}
\int_{Y^*}e_{k\alpha\beta}(x,y)\frac{ \partial \varphi^{mh}}{\partial y_k} s_{\alpha \beta,y}({\bf w}^{ij})~dy=\int_{Y^*} c_{\alpha \beta kl}(x,y) s_{kl,y}(\Sigma^{mh}+{\bf w}^{mh}) s_{\alpha \beta,y}({\bf w}^{ij})~dy.
\end{equation}
Multiplying the second equation by $\varphi^{mh}$, and integrating by parts, we have
\begin{equation}
\int_{Y^*}e_{k\alpha\beta}(x,y)\frac{ \partial \varphi^{mh}}{\partial y_k} s_{\alpha \beta,y}(\Sigma^{ij}+{\bf w}^{ij})~dy=\int_{Y^*} d_{k\alpha}(x,y) \frac{ \partial \varphi^{mh}}{\partial y_k}\frac{ \partial \varphi^{ij}}{\partial y_{\alpha}}~dy.
\end{equation}
Regrouping these results, and using the definition (\ref {Ch}), we derive
\begin{eqnarray}
c_{ijmh}^H & = &\int_{Y^*}c_{ijkl}(x,y)s_{kl,y}(\Sigma^{mh}+ {\bf w}^{mh})~dy+\int_{Y^*}e_{kij}(x,y) \frac{ \partial \varphi^{mh}}{\partial y_k}~dy, \nonumber\\
& = & \int_{Y^*}c_{ijkl}(x,y)s_{kl,y}(\Sigma^{mh}+ {\bf w}^{mh})~dy\nonumber\\
&+&\int_{Y^*} c_{\alpha \beta kl}(x,y) s_{kl,y}(\Sigma^{mh}+{\bf w}^{mh}) s_{\alpha \beta,y}({\bf w}^{ij})~dy + \int_{Y^*} d_{k\alpha}(x,y) \frac{ \partial \varphi^{mh}}{\partial y_k}\frac{ \partial \varphi^{ij}}{\partial y_{\alpha}}~dy,\nonumber\\
& = & \int_{Y^*}c_{\alpha\beta kl}(x,y)s_{kl,y}(\Sigma^{mh}+ {\bf w}^{mh})s_{\alpha \beta,y}(\Sigma^{ij})~dy+
\int_{Y^*} d_{k\alpha}(x,y) \frac{ \partial \varphi^{mh}}{\partial y_k}\frac{ \partial \varphi^{ij}}{\partial y_{\alpha}}~dy\nonumber\\
&+&\int_{Y^*} c_{\alpha \beta kl}(x,y) s_{kl,y}(\Sigma^{mh}+{\bf w}^{mh}) s_{\alpha\beta,y}({\bf w}^{ij})~dy,\nonumber\\
& = & \int_{Y^*}c_{\alpha\beta kl}(x,y) s_{kl,y}(\Sigma^{mh}+ {\bf w}^{mh})s_{\alpha \beta ,y}(\Sigma^{ij}+{\bf w}^{ij})~dy\nonumber\\&+&\int_{Y^*} d_{k\alpha}(x,y) \frac{ \partial \varphi^{mh}}{\partial y_k}\frac{ \partial \varphi^{ij}}{\partial y_{\alpha}}~dy.
\label{symc}
\end{eqnarray}
It is immediate from above, that the coefficients of elasticity tensor satisfies
$$c^H_{ijmh}=c^H_{mhij}$$
This is the end of the proof  of the first section of proposition i.e of symmetry.\\ 
We now study the ellipticity of the coefficients of the elasticity tensor, we recall $c^H_{ijkl}$ is elliptic, if
for all the second order tensor $X_{ij}$ symetric ($X_{ij}=X_{ji}$), we have 
$$\exists ~\alpha^H_c >0, ~~c^H_{ijkl}~X_{ij}~X_{kl} \geq \alpha^H_c~X_{ij}~X_{ij} $$
Consider expression (\ref {c}) of the tensor $c^H_{ijmh}$, we have
\begin{eqnarray*}
c^H_{ijmh}~X_{ij}~X_{mh}& = &\int_{Y^*}c_{ijmh}(x,y)~X_{ij}~X_{mh}~dy+\int_{Y^*}c_{ijkl}(x,y)s_{kl,y}({\bf w}^{mh})~X_{ij}~X_{mh}~dy\\
& +&\int_{Y^*}e_{kij}(x,y)\frac{ \partial \varphi^{mh}}{\partial y_k}~X_{ij}~X_{mh}~dy,\\
&= &\int_{Y^*}c_{ijmh}(x,y)~X_{ij}~X_{mh}~dy+\int_{Y^*}c_{ijkl}(x,y)s_{kl,y}({\bf w}^{mh}X_{mh})X_{ij}~dy\\
&+&\int_{Y^*}e_{kij}(x,y)\frac{ \partial (\varphi^{mh}X_{mh})}{\partial y_k}~X_{ij}~dy,\\
&=&\int_{Y^*}c_{ijkl}(x,y)\Big(s_{kl,y}({\bf w})+P_{kl}\Big)P_{ij}~dy+\int_{Y^*}e_{kij}(x,y)\frac{ \partial \zeta}{\partial y_k}P_{ij}~dy,
\end{eqnarray*}
where ${\bf w}={\bf w}^{mh}X_{mh},~\zeta=\varphi^{mh}X_{mh}$ and $P_{ij}=\tau^{ij}_{mh}X_{mh}=X_{ij}$ . Therfore the couple 
$({\bf w},\zeta)$ is a saddle point of the functional $J$ defined by 
$$
({\bf v}, \psi) \rightarrow J({\bf v},\psi)
$$
$$
J({\bf v},\psi)= \frac{1}{2}\int_{Y^*}\Big\{c_{ijkl}\Big(s_{ij,y}({\bf v})+P_{ij}\Big)\Big(s_{kl,y}({\bf v})+P_{kl}\Big)+2e_{kij}\frac{\partial \psi}{\partial y_k}\Big(s_{ij,y}({\bf v})+P_{ij}\Big)-d_{ij}\frac{ \partial \psi}{\partial y_i}\frac{ \partial \psi}{\partial y_j}\Big\}~dy,
$$
By definition of the saddle point, we have
$$
J({\bf w},\psi) \leq J({\bf w},\zeta) \leq J({\bf v},\zeta) ~~~\forall ({\bf v},\psi)~\mbox{periodic functions.}
$$
Or  
$$ J({\bf w},0)=\frac{1}{2}\int_{Y^*}c_{ijkl}(x,y)\Big(s_{ij,y}({\bf w})+P_{ij}\Big)\Big(s_{kl,y}({\bf w})+P_{kl}\Big)~dy.$$
But, if we use the first equation of system (\ref {pb1}), we obtain
\begin{eqnarray*}
c^H_{ijmh}~X_{ij}~X_{mh}&=&2J({\bf w},\zeta),\\
&\geq& J({\bf w},0),\\
&=&\frac{1}{2}\int_{Y^*}c_{ijkl}(x,y)\Big(s_{ij,y}({\bf w})+P_{ij}\Big)\Big(s_{kl,y}({\bf w})+P_{kl}\Big)~dy, \\
&>&0,
\end{eqnarray*}
We have the homogenized elastcity tensor $\mathcal{C}^H=(c^H_{ijkl})$, which is elliptic.\\
Now we give a results concerning some properties of dielectric homogenized tensor.
\begin{prop}
The coefficients of dielectric homogenized tensor $\mathcal{D}^H=(d^H_{in})$ defined by (\ref{d}) satisfy:
\\\\
a) {$d^H_{in}=d^H_{ni}, \quad \forall 1 \leq i,n \leq 3, $ }
\\
b) There exists $\alpha^H_d >0$, such that for all vector $\xi$, 
$$ d^H_{in}~\xi_i~\xi_n\geq \alpha^H_d~\xi_i~\xi_i$$ 
\end{prop}
{\bf Proof.}\\
By analgy, we transform these coefficients to obtain a symmetric form
$$d^H_{in}=d^H_{ni}$$ 
The problem (\ref {pb2}), can be rewritten as
\begin{eqnarray}
\left\{
\begin{array}{lcllll}
\displaystyle-\frac{\partial}{\partial y_j}\Big\{c_{ijkl}(x,y)s_{kl,y}({\bf q}^{n})+e_{kij}(x,y)\frac{\partial (y_{n}+\psi^{n})}{\partial y_k}\Big\}&=0~\mbox{in}~Y^*,\\ \\
\displaystyle-\frac{\partial}{\partial y_i}\Big\{-e_{ikl}(x,y)s_{kl,y}({\bf q}^{n})+d_{ij}(x,y)\frac{\partial (y_{n}+\psi^{n})}{\partial y_j}\Big\}&=0~\mbox{in}~Y^*,\\ \\
\displaystyle{\bf q}^{n},~~\psi^{n}&Y^*-\mbox{periodics}.
\end{array}
\right.
\label{pb24}
\end{eqnarray}
Introducing $({\bf q}^i,\psi^i)$, in the solution of this local problem
\begin{eqnarray}
\left\{
\begin{array}{lclll}
\displaystyle-\frac{\partial}{\partial y_j}\Big\{c_{ijkl}(x,y)s_{kl,y}({\bf q}^{i})+e_{kij}(x,y)\frac{\partial (y_{i}+\psi^{i})}{\partial y_k}\Big\}&=0~\mbox{in}~Y^*,\\ \\
\displaystyle-\frac{\partial}{\partial y_i}\Big\{-e_{ikl}(x,y)s_{kl,y}({\bf q}^{i})+d_{ij}(x,y)\frac{\partial (y_{i}+\psi^{i})}{\partial y_j}\Big\}&=0 ~\mbox{in}~Y^*,\\ \\
\displaystyle{\bf q}^{i},~~\psi^{i}&Y^*-\mbox{periodics}.
\end{array}
\right.
\label{pb26}
\end{eqnarray}
We can rewrite these coefficients of electric tensor in form given as
\begin{equation}
d_{in}^H  =\int_{Y^*}-e_{ikl}(x,y)s_{kl,y}({\bf q}^{n})~dy+\int_{Y^*}d_{ij}(x,y)\frac {\partial (y_{n}+\psi^n)}{\partial y_j}~dy.
\label{dhh}
\end{equation}
The first term of the second integral of precedent expression, is evaluated as follows
\begin{eqnarray}
\int_{Y^*}-e_{ikl}(x,y)s_{kl,y}({\bf q}^{n})~dy&=&\int_{Y^*}-e_{\alpha kl}(x,y)s_{kl,y}({\bf q}^{n})\delta_{\alpha i}~dy,\nonumber \\
&=& \int_{Y^*}-e_{\alpha kl}(x,y)s_{kl,y}({\bf q}^{n})\frac{\partial y_i}{\partial y_{\alpha}}~dy,\nonumber \\
&=&\int_{Y^*}-e_{\alpha kl}(x,y)s_{kl,y}({\bf q}^{n})\frac{\partial \psi^i}{\partial y_{\alpha}}~dy\nonumber \\
&-&\int_{Y^*}-e_{\alpha kl}(x,y)s_{kl,y}({\bf q}^{n})\frac{\partial (y_i+\psi^i)}{\partial y_{\alpha}}~dy.\nonumber\\
\label{thp}
\end{eqnarray}
Using the variationnal formulation of the second equation of system (\ref {pb24}), and chossing a test function $\varphi=\psi^i$
, we obtain
$$
\int_{Y^*}-e_{\alpha kl}(x,y)s_{kl,y}({\bf q}^{n})\frac{\partial \psi^i}{\partial y_{\alpha}}~dy=\int_{Y^*} d_{\alpha j}(x,y)\frac{\partial (y_n+\psi^n)}{\partial y_j}\frac{\partial \psi^i}{\partial y_{\alpha}}~dy.
$$
Let us now consider the second integral in (\ref {thp}). Multiplying the first equation of system (\ref {pb26}) by
 $\phi^n$, and integrating by parts, we have
$$
\int_{Y^*}-e_{\alpha kl}(x,y)s_{kl,y}({\bf q}^{n})\frac{\partial (y_i+\psi^i)}{\partial y_{\alpha}}~dy=-\int_{Y^*} c_{kl \alpha \beta}(x,y)s_{\alpha \beta, y}({\bf q}^i)s_{kl,y}(\varphi^n)~dy.
$$ 
Finality, we regroup these lasts results, and using the definition (\ref {dhh}), we obtain
\begin{eqnarray}
d^H_{in}&=&\int_{Y^*}-e_{ikl}(x,y)s_{kl,y}({\bf q}^{n})~dy +\int_{Y^*}d_{ij}(x,y)\frac {\partial (y_{n}+\psi^n)}{\partial y_j} ~dy,\nonumber \\
&=&\int_{Y^*} d_{\alpha j}(x,y)\frac{\partial (y_n+\psi^n)}{\partial y_j} \frac{\partial \psi^i}{\partial y_{\alpha}}~dy-\int_{Y^*} c_{kl \alpha \beta}(x,y)s_{\alpha \beta, y}({\bf q}^i)s_{kl,y}({\bf q}^n)~dy\nonumber \\
&+&\int_{Y^*}d_{ij}(x,y)\frac {\partial (y_{n}+\psi^n)}{\partial y_j}~dy,\nonumber \\
&=&\int_{Y^*} d_{\alpha j}(x,y)\frac{\partial (y_n+\psi^n)}{\partial y_j} \frac{\partial \psi^i}{\partial y_{\alpha}}~dy-\int_{Y^*} c_{kl \alpha \beta}(x,y)s_{\alpha \beta, y}({\bf q}^i)s_{kl,y}({\bf q}^n)~dy\nonumber \\
&+&\int_{Y^*}d_{\alpha j}(x,y)\frac {\partial (y_{n}+\psi^n)}{\partial y_j} \frac{\partial y_i}{\partial y_{\alpha}}~dy,\nonumber \\
&=&\int_{Y^*}d_{\alpha j}(x,y)\frac {\partial (y_{n}+\psi^n)}{\partial y_j} \frac{\partial (y_i+\psi^i)}{\partial y_{\alpha}}~dy-\int_{Y^*} c_{kl \alpha \beta}(x,y) s_{\alpha \beta, y}({\bf q}^i)s_{kl,y}({\bf q}^n)~dy.\nonumber \\
\label{hjk}
\end{eqnarray}
It is clear from that the coefficients of electric tensor is symmetric.\\
Now we are interested in the ellipticity of this tensor, recall $d^H_{in}$ is elliptic, if
for all vector $X_i$, we have
$$ \exists ~\alpha^H_d >0, ~~d^H_{in}~X_i~X_n\geq \alpha^H_d~X_i~X_i.$$
We consider the expression (\ref {d}) of tensor $d^H_{in}$, we derive
\begin{eqnarray*}
d^H_{in}~X_i~X_n&=&\int_{Y^*}d_{in}(x,y)~X_i~X_n~dy-\int_{Y^*}e_{ikl}(x,y)s_{kl,y}({\bf q}^{n})~X_i~X_n~dy\\
&+&\int_{Y^*}~d_{ij}(x,y)\frac{\partial \psi^{n}}{\partial y_j}~X_i~X_n~dy, \\
&=&\int_{Y^*}d_{in}(x,y)~X_i~X_n~dy-\int_{Y^*}e_{ikl}(x,y)s_{kl,y}({\bf q}^{n}X_n)~X_i~dy  \\
&+&\int_{Y^*}d_{ij}(x,y)\frac{\partial (\psi^{n}X_n)}{\partial y_j}~X_i~dy,\\
&=&\int_{Y^*}d_{ij}(x,y)\Big(Q_j+\frac{\partial \xi}{\partial y_j}\Big) Q_i~dy-\int_{Y^*}e_{ikl}(x,y)s_{kl,y}(\varsigma)Q_i~dy,
\end{eqnarray*}
where $\xi=\psi^{n}X_n,~\varsigma=\varphi^{n}X_n$ and $Q_i=\delta_{in}X_n=X_i$. Or the couple $(\xi,\varsigma)$ is a saddle point of the functionnal $G$ defined by :
$$
({\bf v}, \psi) \rightarrow G({\bf v},\psi)
$$
$$
G({\bf v},\psi)=\frac{1}{2}\int_{Y^*}\Big\{c_{ijkl}s_{ij,y}({\bf v})s_{kl,y}({\bf v})+2e_{kij}s_{ij,y}({\bf v})\Big(Q_k+\frac{\partial \psi }{\partial y_k}\Big)-d_{ij}\Big(Q_i+\frac{ \partial \psi}{\partial y_i}\Big)\Big(Q_j+\frac{ \partial \psi}{\partial y_j}\Big)\Big\}~dy.
$$
By definition of the saddle point, we have
$$
G(\xi,\psi) \leq G(\xi,\varsigma) \leq G({\bf v},\varsigma) ~~~\forall ({\bf v},\psi)~\mbox{periodic functions.}
$$
Or
$$   
G(0,\varsigma)=-\frac{1}{2}\int_{Y^*}d_{ij}(x,y)(Q_i+\frac{ \partial \varsigma}{\partial y_i})(Q_j+\frac{ \partial \varsigma}{\partial y_j})~dy
$$
But, if we use the second equation of system (\ref {pb2}), we have
\begin{eqnarray*}
G(0,\varsigma)&=&-\frac{1}{2}\int_{Y^*}d_{ij}(x,y)\Big(Q_j+\frac{ \partial \varsigma}{\partial y_j}\Big)Q_i~dy-\frac{1}{2}\int_{Y^*}d_{ij}(x,y)\Big(Q_j+\frac{ \partial \varsigma}{\partial y_j}\Big)\frac{ \partial \varsigma}{\partial y_i}~dy,\\
&=&-\frac{1}{2}\int_{Y^*}d_{ij}(x,y)\Big(Q_j+\frac{ \partial \varsigma}{\partial y_j}\Big)Q_i~dy+\frac{1}{2}\int_{Y^*}e_{ikl}(x,y)s_{kl,y}(\varsigma)~X_i~dy,\\
&=&-\frac{1}{2}d^H_{in}~X_i~X_n,\\
&<&0.
\end{eqnarray*}
We have the dielectric homoginized tensor which is elliptic.\\
Now we give the results concerning some properties of piezoelectric homogenized tensor.
\begin{prop}
The coefficients of piezoelectric homogenized tensor $\mathcal{E}^H=(e^H_{nij})$ defined by (\ref{e}) satisfy:
\\
$$ e^H_{nij}=e^H_{nji}$$
Moreover, we have the identity 
$$
e^H_{nij}=f^H_{nij}
$$
\end{prop}
{\bf Proof.} \\
By definition of coefficients $e^H_{nij}$ (using the fact that $c_{ijkl}(x,y)=c_{jikl}(x,y)$, $e_{kji}(x,y)=e_{kij}(x,y)$), 
we have
\begin{eqnarray}
e_{nij}^H & =& \int_{Y^*} \Big \{c_{ijkl}(x,y)s_{kl,y}({\bf q}^{n})+e_{kij}(x,y)\Big (\delta_{kn}+\frac{\partial \psi^{n}}{\partial y_k}\Big )\Big \}~dy, \nonumber \\
& =& \int_{Y^*} \Big \{c_{jikl}(x,y)s_{kl,y}({\bf q}^{n})+e_{kji}(x,y)\Big (\delta_{kn}+\frac{\partial \psi^{n}}{\partial y_k}\Big )\Big \}~dy, \nonumber \\
& =& e^H_{nji}.
\end{eqnarray}
We can rewrite the coefficients $e^H_{nij}$, as given as under 
\begin{eqnarray}
e_{nij}^H & =& \int_{Y^*} \Big \{c_{ijkl}(x,y)s_{kl,y}({\bf q}^{n})+e_{kij}(x,y)\Big (\delta_{kn}+\frac{\partial \psi^{n}}{\partial y_k}\Big )\Big \}~dy, \nonumber \\
& =& \int_{Y^*} \Big \{e_{nij}(x,y)+e_{kij}(x,y)\frac{\partial \psi^{n}}{\partial y_k}+c_{ijkl}(x,y)s_{kl,y}({\bf q}^{n})\Big \}~dy. 
\label{ef}
\end{eqnarray}
Same, we can rewrite the coefficients $f^H_{nij}$, as form givens under
\begin{eqnarray}
f_{nij}^H & =& \int_{Y^*}\Big \{e_{nkl}(x,y)\Big (\tau^{kl}_{ij}+s_{kl,y}({\bf w}^{ij})\Big )-d_{nt}(x,y)\frac{\partial \varphi^{ij}}{\partial y_t}\Big \}~dy, \nonumber\\
&=&\int_{Y^*}\Big \{e_{nij}(x,y)+e_{nkl}(x,y)s_{kl,y}({\bf w}^{ij})-d_{nt}(x,y)\frac{\partial \varphi^{ij}}{\partial y_t}\Big \}~dy.
\end{eqnarray}
Using the two variationnals formulations corresponding of problems (\ref {pb1}) and (\ref {pb2}), and chossing the appropriete test functions, we can directly prove
as 
$e^H_{nij}=f^H_{nij}$.
\\
Finallity, using the three last propositions, we can purpose the altarnative form of the principal convergence theorem \\
{\bf Theorem-Bis ({\it the altarnative form})}\\
{\it Set $({\bf u}, \varphi)$ solution of the two-scale homogenized problem (\ref {pbhom1})-(\ref {pbhom2})-(\ref {pbhom3}), then $({\bf u}, \varphi)$ is defined by that the solution of this homogenized problem
\begin{eqnarray}
\left\{
\begin{array}{lcll}
\displaystyle -{\bf div}~\sigma^H({\bf u},\varphi)& =& \theta~{\bf f} &\mbox{in}~\Omega,\\ \\
\displaystyle -{\bf div}~D^H({\bf u},\varphi)& = &0 &\mbox{in}~\Omega,
\end{array}
\right.
\label{pb3}
\end{eqnarray}
where the boundary conditions 
\begin{equation}
\left\{
\begin{array}{lcl}
{\bf u}(x)& = ~{\bf 0}  &\mbox{on}~\partial\Omega,\\ \\ 
\varphi(x)&= ~0  &\mbox{on}~\partial\Omega, 
\end{array}
\right.
\end{equation}
$\sigma^H_{ij}$ and $D^H_{i}$ are defined by the homogenized constitutive law
\begin{equation}
\left\{
\begin{array}{lcl}
\displaystyle \sigma^H_{ij}({\bf u},\varphi)&=& \displaystyle c_{ijmh}^H s_{mh,x}({\bf u})+ e_{nij}^H \frac{\partial \varphi}{\partial x_{n}},\\ \\
\displaystyle D^H_{i}({\bf u},\varphi) &= &\displaystyle-e_{imh}^H s_{mh,x}({\bf u})+d_{in}^H \frac{\partial \varphi}{\partial x_{n}},
\end{array}
\right.
\end{equation}
the homogenized coefficients $c^H_{ijkl}$, $e^H_{nij}$ and $d^H_{ij}$ are defined respectively by (\ref {c}), (\ref {e}) and (\ref {d}).}
\section{Correctors result}
The corrector results are obtained easily by the two-scale convergence method. The objective of the next theorem justify rigorously  the two first terms in the usual asymptotic expansion of the solution.
\\
Following the idea of Allaire \cite {All2}, we introduce the following definition
\begin{defi}
We call $ \psi(x,y)$ an {\it admissible test function}, if it is $Y$-periodic, and satisfies the following relation
\begin{equation}
\lim_{\varepsilon \rightarrow 0}\int_{\Omega}\psi(x,\frac{x}{\varepsilon})^2 dx = \int_{\Omega}\int_{Y} \psi(x,y)^2 dx~dy
\label{khj}
\end{equation}
\label{DefAl}
\end{defi}
Here we recall the Allaire's lemma\\
\begin{lem} (Allaire \cite {All2}) :\\
Let the function $\psi(x,y) \in L^2(\Omega;C_{\sharp}(Y))$, then $\psi(x,y)$ is an admissible test function in the sense of Definition \ref{DefAl}.
\label{LemAll}
\end{lem}
Using this lemma, we obtain the following proposition.
\begin{prop}
The two functions $s_{ij,y}({\bf u}_1(x,y))$ and $\partial_{i,y} \varphi_1(x,y)$ are admissible test functions in the sense of Definition \ref{DefAl}. 
\end{prop}
{\bf Proof.}\\
By definition we have
$$
{\bf u}_1(x,y) = s_{mh,x}({\bf u}(x)){\bf w}^{mh}(y)+\frac{\partial \varphi(x)}{\partial x_{n}} {\bf q}^{n}(y),
$$
$$
\varphi_1(x,y) = s_{mh,x}({\bf u}(x)) \varphi^{mh}(y)+\frac{\partial \varphi(x)}{\partial x_{n}} \psi^{n}(y),
$$
we obtain
$$
s_{ij,y}({\bf u}_1(x,y)) = s_{mh,x}({\bf u}(x))s_{ij,y}({\bf w}^{mh}(y))+\frac{\partial \varphi(x)}{\partial x_{n}} s_{ij,y}({\bf q}^{n}(y)),
$$
$$
\partial_{i,y} \varphi_1(x,y) = s_{mh,x}({\bf u}(x)) \partial_{i,y} \varphi^{mh}(y)+\frac{\partial \varphi(x)}{\partial x_{n}} \partial_{i,y} \psi^{n}(y).
$$
Using Lemma \ref{LemAll}, $s_{ij,y}({\bf u}_1(x,y))$ and $\partial_{i,y} \varphi_1(x,y)$ are the admissible test functions in the sense of Definition \ref{DefAl}. 
\begin{theo}
We have these two strong convergence results
$$
\left\{
\begin{array}{lcl}
\displaystyle\stackrel{\sim}{{\bf u}^{\varepsilon}}(x)-\chi(\frac{x}{\varepsilon})[{\bf u}(x)+{\bf u}_1(x,\frac{x}{\varepsilon})]\rightarrow 0 &\mbox{strongly in }~{\bf H}^1_0(\Omega)\\
\displaystyle\stackrel{\sim}{\varphi^{\varepsilon}}(x)-\chi(\frac{x}{\varepsilon})[\varphi(x)+ \varphi_1(x,\frac{x}{\varepsilon})]\rightarrow 0 &\mbox{strongly in }~H^1_0(\Omega)
\end{array}
\right.
$$
\end{theo}
{\bf Proof.} We consider the variational formulation under the following form 
\begin{equation}
\int_{\Omega_{\varepsilon}} \{[c^{\varepsilon}_{ijkl} s_{kl}({\bf u}^{\varepsilon}) +e^{\varepsilon}_{kij} \partial_{k}\varphi^{\varepsilon}]s_{ij}({\bf v})+[-e^{\varepsilon}_{ikl} s_{kl}({\bf u}^{\varepsilon})+d^{\varepsilon}_{ij}\partial_{j} \varphi^{\varepsilon}]\partial_{i} \psi \}~dx = \int_{\Omega_{\varepsilon}} f_i(x) v_i(x)~dx
\label{thg}
\end{equation}
Chosing ${\bf v}={\bf u}^{\varepsilon} \mbox{ and }\psi=\varphi^{\varepsilon}$ in (\ref {thg}), we obtain
$$
\int_{\Omega_{\varepsilon}} \{[c^{\varepsilon}_{ijkl} s_{kl}({\bf u}^{\varepsilon}) +e^{\varepsilon}_{kij} \partial_{k} \varphi^{\varepsilon}]s_{ij}({\bf u}^{\varepsilon})+[-e^{\varepsilon}_{ikl} s_{kl}({\bf u}^{\varepsilon})+d_{ij}^{\varepsilon}\partial_{j} \varphi^{\varepsilon}]\partial_{i} \varphi^{\varepsilon} \}~dx = \int_{\Omega_{\varepsilon}} f_i(x) u_i^{\varepsilon}(x) ~dx
$$
By simplification, we get
\begin{equation}
\int_{\Omega_{\varepsilon}} \{ c^{\varepsilon}_{ijkl}(x)s_{ij}({\bf u}^{\varepsilon})(x)s_{kl}({\bf u}^{\varepsilon})(x)+d_{ij}^{\varepsilon}(x) \partial_i \varphi^{\varepsilon}(x)\partial_j \varphi^{\varepsilon}(x)\}~dx = \int_{\Omega_{\varepsilon}} f_i(x) u_i^{\varepsilon}(x) ~dx
\label{er}
\end{equation}
By applying (\ref {er}), we can write 
\begin{eqnarray*}
&\displaystyle\int_{\Omega}&c^{\varepsilon}_{ijkl}\Big \{\stackrel{\sim}{s}_{ij}({\bf u}^{\varepsilon})-\chi(\frac{x}{\varepsilon})\Big[s_{ij,x}({\bf u})+s_{ij,y}({\bf u}_1)\Big]\Big \}\Big \{\stackrel{\sim}{s}_{kl}({\bf u}^{\varepsilon})-\chi(\frac{x}{\varepsilon})\Big[s_{kl,x}({\bf u})+s_{kl,y}({\bf u}_1)\Big]\Big \}dx\nonumber \\\\
&+\displaystyle \int_{\Omega}&d^{\varepsilon}_{ij}(x)\Big \{\stackrel{\sim}{\partial_i}\varphi^{\varepsilon}(x)-\chi(\frac{x}{\varepsilon})\Big[\partial_{i,x} \varphi+\partial_{i,y} \varphi_1\Big]\Big\}\Big\{ \stackrel{\sim}{\partial_j}\varphi^{\varepsilon}(x)-\chi(\frac{x}{\varepsilon})\Big[\partial_{j,x} \varphi+\partial_{j,y} \varphi_1\Big]\Big\}dx\nonumber\\ \\
&=&\int_{\Omega}f_i(x)\stackrel{\sim}{ u_i^{\varepsilon}}(x)~dx\nonumber\\\\
&+&\int_{\Omega} c^{\varepsilon}_{ijkl}(x)\chi(\frac{x}{\varepsilon})\Big [s_{ij,x}({\bf u}(x))+s_{ij,y}({\bf u}_1(x,\frac{x}{\varepsilon}))\Big ]\Big [s_{kl,x}({\bf u}(x))+s_{kl,y}({\bf u}_1(x,\frac{x}{\varepsilon}))\Big ]~dx\nonumber \\\\
&+&\int_{\Omega} d^{\varepsilon}_{ij}(x)\chi(\frac{x}{\varepsilon})\Big [\partial_{i,x}\varphi(x)+\partial_{i,y} \varphi_1(x,\frac{x}{\varepsilon})\Big ]\Big [\partial_{j,x}\varphi(x)+\partial_{j,y} \varphi_1(x,\frac{x}{\varepsilon})\Big ]~dx\nonumber \\\\
&-&2\int_{\Omega}c^{\varepsilon}_{ijkl}(x)\chi(\frac{x}{\varepsilon})\stackrel{\sim}{s_{ij}}({\bf u}^{\varepsilon})\Big [s_{kl,x}({\bf u}(x))+s_{kl,y}({\bf u}_1(x,\frac{x}{\varepsilon}))\Big ]~dx\nonumber \\\\
&-&2\int_{\Omega}d^{\varepsilon}_{ij}(x)\chi(\frac{x}{\varepsilon})\stackrel{\sim}{\partial_i}\varphi^{\varepsilon}(x)\Big [\partial_{j,x} \varphi(x)+\partial_{j,y} \varphi_1(x,\frac{x}{\varepsilon})\Big ]~dx.
\end{eqnarray*}
Using the ellipticity property of the elastic $(c^{\varepsilon}_{ijkl})$ and dielectric $(d^{\varepsilon}_{ij})$ tensors, we get
\begin{eqnarray*}
&\alpha_c& \parallel \stackrel{\sim}{s}_{ij}({\bf u}^{\varepsilon})-\chi(\frac{x}{\varepsilon})s_{ij,x}({\bf u}(x))-\chi(\frac{x}{\varepsilon})s_{ij,y}({\bf u}_1(x,\frac{x}{\varepsilon}))\parallel^2_{{\bf L}^2(\Omega)}\\\\
+&\alpha_d& \parallel\stackrel{\sim}{\partial_i}\varphi^{\varepsilon}(x)-\chi(\frac{x}{\varepsilon})\partial_{i,x} \varphi(x)-\chi(\frac{x}{\varepsilon})\partial_{i,y} \varphi_1(x,\frac{x}{\varepsilon}) \parallel^2_{L^2(\Omega)}\\\\
&\leq& \int_{\Omega} f_i(x)\stackrel{\sim}{u^{\varepsilon}_i}(x)~dx\\\\
&+&\int_{\Omega}c^{\varepsilon}_{ijkl}(x)\chi(\frac{x}{\varepsilon})\Big [s_{ij,x}({\bf u}(x))+s_{ij,y}({\bf u}_1(x,\frac{x}{\varepsilon}))\Big ]\Big [s_{kl,x}({\bf u}(x))+s_{kl,y}({\bf u}_1(x,\frac{x}{\varepsilon}))\Big ]~dx\\\\
&+&\int_{\Omega}d^{\varepsilon}_{ij}(x)\chi(\frac{x}{\varepsilon})\Big [\partial_{i,x}\varphi(x)+\partial_{i,y} \varphi_1(x,\frac{x}{\varepsilon})\Big ]\Big [\partial_{j,x}\varphi(x)+\partial_{j,y} \varphi_1(x,\frac{x}{\varepsilon})\Big ]~dx\\\\
&-2&\int_{\Omega}c^{\varepsilon}_{ijkl}(x)\chi(\frac{x}{\varepsilon})\stackrel{\sim}{s_{ij}}({\bf u}^{\varepsilon})\Big [s_{kl,x}({\bf u}(x))+s_{kl,y}({\bf u}_1(x,\frac{x}{\varepsilon}))\Big ]~dx\\\\
&-2&\int_{\Omega}d^{\varepsilon}_{ij}(x)\chi(\frac{x}{\varepsilon})\stackrel{\sim}{\partial_i}\varphi^{\varepsilon}(x)\Big [\partial_{j,x} \varphi(x)+\partial_{j,y} \varphi_1(x,\frac{x}{\varepsilon})\Big ]~dx.
\end{eqnarray*}
Using the fact that $s_{ij,y}({\bf u}_1(x,y))$ and $\partial_{i,y} \varphi_1(x,y)$ are the admissible test functions and taking the limit in the sense of the two-scale convergence in the second right-hand side of the inequality, we obtain
\begin{eqnarray}
 &\displaystyle\alpha_c\lim_{\varepsilon \rightarrow 0} \parallel \stackrel{\sim}{s_{ij}}({\bf u}^{\varepsilon})-\chi(\frac{x}{\varepsilon})\{s_{ij,x}({\bf u}(x))-s_{ij,y}({\bf u}_1(x,\frac{x}{\varepsilon}))\} \parallel^2_{{\bf L}^2(\Omega)}&\nonumber\\
&+\displaystyle \alpha_d \lim_{\varepsilon \rightarrow 0}\parallel\stackrel{\sim}{\partial_i}\varphi^{\varepsilon}(x)-\chi(\frac{x}{\varepsilon})\{\partial_{i,x} \varphi(x)-\partial_{i,y} \varphi_1(x,\frac{x}{\varepsilon})\} \parallel^2_{L^2(\Omega)} &\nonumber\\
&\displaystyle \leq \int_{\Omega}\int_{Y^*} f_i(x) u_i(x)~dx~dy &\nonumber\\
&-\displaystyle \int_{\Omega}\int_{Y^*}c_{ijkl}(x,y)[s_{ij,x}({\bf u}(x))+s_{ij,y}({\bf u}_1(x,y))][s_{kl,x}({\bf u}(x))+s_{kl,y}({\bf u}_1(x,y))]~dx~dy&\nonumber\\
&-\displaystyle \int_{\Omega}\int_{Y^*} d_{ij}(x,y)[\partial_{i,x}\varphi(x)+\partial_{i,y}\varphi_1(x,y)][\partial_{j,x}\varphi(x)+\partial_{j,y}\varphi_1(x,y)]~dx~dy
\label{tat}
\end{eqnarray}
Recalling the form of the two-scale homogenized problem (\ref {pbhom1})-(\ref {pbhom2})-(\ref {pbhom3}), we observe that the right-hand side of the inequality (\ref {tat}) vanishes, so tha, we obtain
$$
\displaystyle\lim_{\varepsilon \rightarrow 0} \parallel \stackrel{\sim}{s_{ij}}({\bf u}^{\varepsilon})-\chi(\frac{x}{\varepsilon})\{s_{ij,x}({\bf u}(x))-s_{ij,y}({\bf u}_1(x,\frac{x}{\varepsilon}))\} \parallel_{{\bf L}^2(\Omega)} =~0
$$
and
$$
\displaystyle\lim_{\varepsilon \rightarrow 0}\parallel\stackrel{\sim}{\partial_i}\varphi^{\varepsilon}(x)-\chi(\frac{x}{\varepsilon})\{\partial_{i,x} \varphi(x)-\partial_{i,y} \varphi_1(x,\frac{x}{\varepsilon})\} \parallel _{L^2(\Omega)}=~0.
$$
\section{Conclusion}
In this work, we have given the new convergence results, and the explicite forms of the elastic, piezoelectric and dielectric homogenized coefficients. The two-scale convergence is applied to our problem yields the strong convergence result on the correctors. This technique of two-scale convergence can handle also other homogenization
problems, in medium which has periodic structure for example the laminated piezocomposite materials or fiber materials (see \cite {Cast} \cite {Feng} \cite{Mech1} \cite{Ruan} \cite{Mech2}). Numerical implementation for perforated, laminated and fiber structures, will be presented in forthcoming publications (see \cite{Mech2}).\\\\ 
{\bf Acknowledgment.} This work has been supported in part by the Ministry for higher education and scientific research of Algeria (University of Oran, Departement of Mathematics). The author is grateful to Professor Bernadette Miara for helpful discussions.


\begin{thebibliography}{16} \frenchspacing
\bibitem{All2} Allaire G. {\it Homogenization and two scale-convergence}, SIAM J. Math. Anal., {\bf 23} (26), (1992) 1482-1518.

\bibitem{All3} Allaire G., Murat F. {\it Homogenization of Neumann problem with non-isolated holes}, Asymptotic. Anal. {\bf 7}, (1993) 81-95.
\bibitem{Berg} Berger H., Gabbet U., K\"{o}ppe H., Rodriguez-Ramos R., Bravo-Castillero J.,  Guinovart-Diaz R., Otero J.A., Maugin G.A.  \emph{ Finite element and asymptotic homogenization methods applied to smart composite materials}. Comp. Mech. {\bf 33}, (2003) 61-67.
\bibitem{Ben} Bensoussan  A., Lions J.L., Papanicolaou G. {\it Asymptotic Analysis for Periodic Structures}, North  Holland, Amsterdam (1978).
\bibitem{Cast} Castillero J.B., Otero J.A., Ramos R.R., Bourgeat A. \emph{ Asymptotic homogenization of laminated piezocomposite materials}, Int. J. Solids Structures. {\bf 35} (1998) 527-541.
\bibitem{Cior} Cioranescu D., Damlamian A., Griso G. {\it Periodic unfolding and homogenization}, C. R. Acad. Sci. Paris, Ser. I {\bf 334} (2002) 99-104.
\bibitem{Gher1} Ghergu M., Griso G., Mechkour H., Miara B. {\it Homog\'en\'eisation de coques minces pi\'ezo\'electriques perfor\'ees}. C. R. Acad. Sci. Paris, Ser. II : M\'ecanique {\bf 333} (2005) 249-255.
\bibitem{Gher3} Ghergu M., Griso G., Labat B., Mechkour H., Miara B., Rohan E., Zidi M. {\it Homog\'en\'eisation et pi\'ezo\'electricit\'e. Aide \`a la conception d'un bio-mat\'eriau}. Annals of University of Craoiva. Math. Comp. Sci. Ser. {\bf 32} (2005): 9-15. 
%\bibitem{Gher2} Ghergu M., Griso G., Mechkour H., Miara B. {\it Homogenization of thin piezoelectric shells}. (2007) To appear.
\bibitem{Feng} Feng M.L., Wu C.C. {\it A study of three-dimensional four-step braided piezo-ceramic composites by the homogenization method}, Comp Scien Tech, {\bf 61} (2001) 1889-1898.
\bibitem{Ngue} Nguetseng G. {\it A general convergence result for a functionnal related to the theory of homogenization}, SIAM J. Math. Anal., {\bf 20}(3), (1989) 608-623.
\bibitem{Mech1} Mechkour H., Miara B. \emph{ Modelling and control of piezoelectric perforated structures}, Proceedings of The Third World Conference On Structural Control. John Wiley, Chichester. F. Casciati : Editor. Vol 3, (2003) 329-336.
\bibitem{Mech2} Mechkour H {\it Homog\'en\'eisation et simulation num\'erique de structures pi\'ezo\'electriques perfor\'ees et lamin\'ees}. PhD thesis, University of Marne-La-Vall\'ee 2004 (in french).
\bibitem{Olei} Oleinik O.A., Shamaev G.A., Yosifian G.A. {\it Mathematical problems in elasticity and homogenization}, North  Holland, Amsterdam (1992).
\bibitem{Ruan} Ruan X., Safari A., Chou T.W. {\it Effective elastic, piezoelectric and dielectric properties of braided fabric composites}, Comp Part A {\bf 30} (1999)1435-1444.
\bibitem{Past} Pastor J \emph{ Homogenization of linear piezoelectric media}. Mech. Res. Comm. {\bf 24}(2), (1997) 145-150.
\bibitem{Tel1} Telega J.J. \emph{ Piezoelectricity and homogenization. Application to biomechanic}, In: Maugin, G.A.(Ed.), Continum Models and Discrete Systems, Vol. 2. Longmam, Harlow, Essex, (1991) 220-229.
\end{thebibliography}
\end{document}